\documentclass[a4paper,10pt]{article}
\usepackage{latexsym}
\usepackage{amssymb}
\usepackage{amsfonts}
\usepackage{amsmath}
\usepackage{indentfirst}
\usepackage{graphicx}
\usepackage{epsfig}
\usepackage[enableskew]{youngtab}
\usepackage{tikz,ifthen}

\newtheorem{prop}{Proposition}[section]
\newtheorem{teor}{Theorem}[section]
\newtheorem{lemma}{Lemma}[section]
\newtheorem{cor}{Corollary}[section]
\newcommand{\nkinN}{n,k\in \mathbf{N}}
\newcommand{\ninN}{n\in \mathbf{N}}

\newcommand{\cvd}{\hfill $\blacksquare$\bigskip}

\newcounter{indice}

\newcommand{\shape}[1]{
\setcounter{indice}{0};
\foreach \i in {#1} {
\addtocounter{indice}{1};
\foreach \x in {1,...,\i} {
\draw[thin] (\x-1,-\theindice+1) rectangle (\x,-\theindice);
}
}
}

\newcommand{\filledshape}[1]{
\setcounter{indice}{0};
\foreach \i in {#1} {
\addtocounter{indice}{1};
\foreach \x in {1,...,\i} {
\filldraw[fill=yellow, draw=black] (\x-1,-\theindice+1) rectangle (\x,-\theindice);
}
}
}

\date{}

\author{Luca Ferrari\thanks{Dipartimento di Sistemi e Informatica, viale Morgagni 65, 50134 Firenze, Italy
{\tt ferrari@dsi.unifi.it}}}

\title{Unimodality and Dyck paths}

\begin{document}

\maketitle

\begin{abstract}
We propose an original approach to the problem of rank-unimodality
for Dyck lattices. It is based on a well known recursive
construction of Dyck paths originally developed in the context of
the ECO methodology, which provides a partition of Dyck lattices
into saturated chains. Even if we are not able to prove that Dyck
lattices are rank-unimodal, we describe a family of polynomials
(which constitutes a  polynomial analog of ballot numbers) and a
succession rule which appear to be useful in addressing such a
problem. At the end of the paper, we also propose and begin a
systematic investigation of the problem of unimodality of
succession rules.
\end{abstract}

\section{Introduction}

In enumerative combinatorics it often happens to discover integer
sequences which are unimodal. A finite sequence $a_0 ,a_1 ,\ldots
a_n$ is said to be \emph{unimodal} when there exists an index
$0\leq i\leq n$ such that $a_0 \leq a_1 \leq \cdots \leq
a_{i-1}\leq a_i \geq a_{i+1}\geq \cdots \geq a_{n-1}\geq a_n$.
Proving that a sequence is unimodal is often a very hard task. A
few papers illustrating some general techniques to tackle this
problem has been published, such as \cite{B,Stanl2} (which provide
a very rich account of several methods to prove unimodality). A
related property is log-concavity. A finite sequence of integers
$a_0 ,a_1 ,\ldots a_n$ is said to be \emph{log-concave} whenever
$a_i ^2 \geq a_{i-1}a_{i+1}$ for all $1\leq i\leq n$. It is not
too difficult to prove that a nonnegative log-concave sequence
having no internal zeroes is unimodal. Since proving log-concavity
is usually easier than proving unimodality, the above result is
often used.

\bigskip

In this paper we will consider unimodality in the context of a
particularly interesting and important combinatorial structure,
i.e. lattice paths. Another paper dealing with unimodality and
lattice paths is \cite{Sa1}, but it seems not to be related to
what we are studying here. For our purposes, a \emph{lattice path}
is a path in the discrete plane starting at the origin of a fixed
Cartesian coordinate system, ending somewhere on the $x$-axis,
never going below the $x$-axis and using only a prescribed set of
steps. This definition is extremely restrictive if compared to
what is called a lattice path in the literature, but it will be
enough for our purposes.

Some very well studied classes of lattice paths are the following:
\begin{itemize}
\item \emph{Dyck paths}, i.e. lattice paths using only steps of
the type $u=(1,1)$ and $d=(1,-1)$;

\item \emph{Motzkin paths}, i.e. lattice paths using only steps of
the type $u=(1,1)$, $h=(1,0)$ and $d=(1,-1)$;

\item \emph{Schr\"oder paths}, i.e. lattice paths using only steps
of the type $u=(1,1)$, $H=(2,0)$ and $d=(1,-1)$;
\end{itemize}

Given a lattice path $P$, the \emph{area} of $P$ is defined to be
the area of the region included between the path $P$ and the
$x$-axis (see figure \ref{area}). So we can consider, for
instance, the distribution of the parameter ``area" over all Dyck
paths of a given length. The following conjecture is not at all a
new one.

\bigskip

\begin{figure}[ht]
\begin{center}
\begin{tikzpicture}
\begin{scope}[scale=0.5]
\draw (0,-4) node {$\bullet$};
\draw (1,-3) node {$\bullet$};
\draw (2,-2) node {$\bullet$};
\draw (3,-3) node {$\bullet$};
\draw (4,-2) node {$\bullet$};
\draw (5,-1) node {$\bullet$};
\draw (6,-2) node {$\bullet$};
\draw (7,-3) node {$\bullet$};
\draw (8,-4) node {$\bullet$};
\draw (9,-3) node {$\bullet$};
\draw (10,-4) node {$\bullet$};
\draw (11,-3) node {$\bullet$};
\draw (12,-2) node {$\bullet$};
\draw (13,-3) node {$\bullet$};
\draw (14,-2) node {$\bullet$};
\draw (15,-3) node {$\bullet$};
\draw (16,-4) node {$\bullet$};

\draw (0,-4) -- (1,-3);
\draw (1,-3) -- (2,-2);
\draw (2,-2) -- (3,-3);
\draw (3,-3) -- (4,-2);
\draw (4,-2) -- (5,-1);
\draw (5,-1) -- (6,-2);
\draw (6,-2) -- (7,-3);
\draw (7,-3) -- (8,-4);
\draw (8,-4) -- (9,-3);
\draw (9,-3) -- (10,-4);
\draw (10,-4) -- (11,-3);
\draw (11,-3) -- (12,-2);
\draw (12,-2) -- (13,-3);
\draw (13,-3) -- (14,-2);
\draw (14,-2) -- (15,-3);
\draw (15,-3) -- (16,-4);

\draw[thin] (0,-4) -- (16,-4);
\draw[thin] (1,-3) -- (7,-3);
\draw[thin] (3,-3) -- (4,-4);
\draw[thin] (4,-2) -- (6,-4);
\draw[thin] (4,-2) -- (6,-2);
\draw[thin] (11,-3) -- (15,-3);
\draw[thin] (1,-3) -- (2,-4);
\draw[thin] (2,-4) -- (3,-3);
\draw[thin] (4,-4) -- (6,-2);
\draw[thin] (6,-4) -- (7,-3);
\draw[thin] (11,-3) -- (12,-4);
\draw[thin] (12,-4) -- (13,-3);
\draw[thin] (13,-3) -- (14,-4);
\draw[thin] (14,-4) -- (15,-3);
\end{scope}
\end{tikzpicture}
\end{center}
\caption{A Dyck path of semilength 8 and area 20. The triangles in the figure
have unit area.\label{area}}
\end{figure}

\bigskip

\textbf{Conjecture.} The sequence $(a_k ^{(n)})_k$ of the number
of Dyck paths of semilength $n$ having area $n+2k$ is unimodal,
for all $n$.

\bigskip

The first few lines of the matrix of the $a_k ^{(n)}$'s are the
following:

\begin{displaymath}
{\scriptsize 
\begin{array}{cccccccccccccccccccccc}
  1 & 0 & 0 & 0 & 0 & 0 & 0 & 0 & 0 & 0 & 0 & 0 & 0 & 0 & 0 & 0 & 0 & 0 & 0 & 0 & 0 & 0 \\
  1 & 1 & 0 & 0 & 0 & 0 & 0 & 0 & 0 & 0 & 0 & 0 & 0 & 0 & 0 & 0 & 0 & 0 & 0 & 0 & 0 & 0 \\
  1 & 2 & 1 & 1 & 0 & 0 & 0 & 0 & 0 & 0 & 0 & 0 & 0 & 0 & 0 & 0 & 0 & 0 & 0 & 0 & 0 & 0 \\
  1 & 3 & 3 & 3 & 2 & 1 & 1 & 0 & 0 & 0 & 0 & 0 & 0 & 0 & 0 & 0 & 0 & 0 & 0 & 0 & 0 & 0 \\
  1 & 4 & 6 & 7 & 7 & 5 & 5 & 3 & 2 & 1 & 1 & 0 & 0 & 0 & 0 & 0 & 0 & 0 & 0 & 0 & 0 & 0 \\
  1 & 5 & 10 & 14 & 17 & 16 & 16 & 14 & 11 & 9 & 7 & 5 & 3 & 2 & 1 & 1 & 0 & 0 & 0 & 0 & 0 & 0 \\
  1 & 6 & 15 & 25 & 35 & 40 & 43 & 44 & 40 & 37 & 32 & 28 & 22 & 18 & 13 & 11 & 7 & 5 & 3 & 2 & 1 & 1 \\
\end{array}
}
\end{displaymath}

Observe that, for $n\geq 3$, none of the sequences $(a_k
^{(n)})_k$ is log-concave. This is trivially seen by observing
that the last three nonzero terms of each sequence are 2,1,1.
However, for specific values of $n$, there are also other terms
that prevent $(a_k ^{(n)})_k$ from being log-concave (as it is
readily seen by inspecting the matrix displayed above). This fact
makes even more intriguing (and surely more difficult) the above
mentioned conjecture.

\bigskip

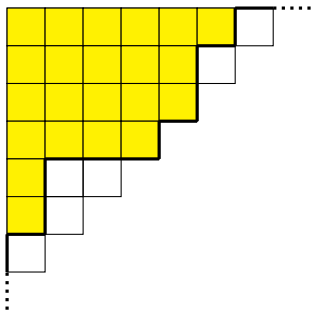
\begin{figure}[ht]
\begin{center}
\begin{tikzpicture}
\begin{scope}[scale=0.5]
\shape{7,6,5,4,3,2,1}
\filledshape{6,5,5,4,1,1}

\draw[very thick,dotted] (0,-8) -- (0,-7);
\draw[very thick] (0,-7) -- (0,-6);
\draw[very thick] (0,-6) -- (1,-6);
\draw[very thick] (1,-6) -- (1,-4);
\draw[very thick] (1,-4) -- (4,-4);
\draw[very thick] (4,-4) -- (4,-3);
\draw[very thick] (4,-3) -- (5,-3);
\draw[very thick] (5,-3) -- (5,-1);
\draw[very thick] (5,-1) -- (6,-1);
\draw[very thick] (6,-1) -- (6,0);
\draw[very thick] (6,0) -- (7,0);
\draw[very thick,dotted] (7,0) -- (8,0);
\end{scope}
\end{tikzpicture}
\end{center}
\caption{The Young diagram corresponding to the Dyck path in figure \ref{area}.\label{Young}}
\end{figure}

This conjecture has first appeared in \cite{Stant}, where it is
however stated in a different language. Indeed, there is a
bijection between Dyck paths of semilength $n$ and Young diagrams
fitting inside the \emph{staircase shape} $(n-1,n-2,\ldots 2,1)$
(see figure \ref{Young}). Moreover, this bijection maps the area of
a Dyck path into the difference between the total area of the
staircase shape and the area of the Young diagram associated with
the path. Thus the above conjecture is formulated by Stanton in
the following form.

\bigskip

\textbf{(Equivalent) Conjecture.} The sequence $(a_k ^{(n)})_k$ of
the number of Young diagrams fitting inside the staircase shape
$(n-1,n-2,\ldots 2,1)$ having area $k$ is unimodal, for all $n$.

\bigskip

The ``path version" of the conjecture of Stanton is due to Bonin,
Shapiro and Simion, who stated it in \cite{BSS}, together with an
analogous conjecture for Schr\"oder paths.

There is still another way to express this unimodality conjecture,
which involves lattices.

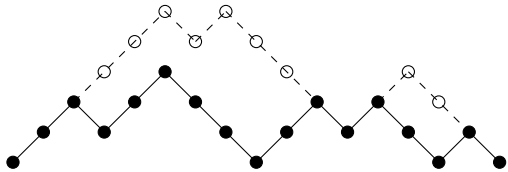
\begin{figure}[ht]
\begin{center}
\begin{tikzpicture}
\begin{scope}[scale=0.4]
\draw (0,0) -- (2,2);
\draw (2,2) -- (3,1);
\draw (3,1) -- (5,3);
\draw (5,3) -- (8,0);
\draw (8,0) -- (10,2);
\draw (10,2) -- (11,1);
\draw (11,1) -- (12,2);
\draw (12,2) -- (14,0);
\draw (14,0) -- (15,1);
\draw (15,1) -- (16,0);
\draw (0,0) [fill] circle (.2);
\draw (1,1) [fill] circle (.2);
\draw (2,2) [fill] circle (.2);
\draw (3,1) [fill] circle (.2);
\draw (4,2) [fill] circle (.2);
\draw (5,3) [fill] circle (.2);
\draw (6,2) [fill] circle (.2);
\draw (7,1) [fill] circle (.2);
\draw (8,0) [fill] circle (.2);
\draw (9,1) [fill] circle (.2);
\draw (10,2) [fill] circle (.2);
\draw (11,1) [fill] circle (.2);
\draw (12,2) [fill] circle (.2);
\draw (13,1) [fill] circle (.2);
\draw (14,0) [fill] circle (.2);
\draw (15,1) [fill] circle (.2);
\draw (16,0) [fill] circle (.2);
\draw[dashed] (2,2) -- (5,5);
\draw[dashed] (5,5) -- (6,4);
\draw[dashed] (6,4) -- (7,5);
\draw[dashed] (7,5) -- (10,2);
\draw[dashed] (12,2) -- (13,3);
\draw[dashed] (13,3) -- (15,1);
\draw (3,3) circle (.2);
\draw (4,4) circle (.2);
\draw (5,5) circle (.2);
\draw (6,4) circle (.2);
\draw (7,5) circle (.2);
\draw (8,4) circle (.2);
\draw (9,3) circle (.2);
\draw (13,3) circle (.2);
\draw (14,2) circle (.2);
\end{scope}
\end{tikzpicture}
\end{center}
\caption{A pair of Dyck paths $P$ (thick) and $Q$ (dashed), with $P<Q$.\label{order}}
\end{figure}

Given a class of paths $\mathcal{P}$, the set $\mathcal{P}_n$ of
all paths in $\mathcal{P}$ having length $n$ can be naturally
endowed with a poset structure, by imposing that $P\leq Q$
whenever $P$ lies weakly below $Q$ (weakly meaning that $P$ and
$Q$ are allowed to have some points in common). See figure
\ref{order} for a ``Dyck" example. It turns out that, for many
interesting classes of paths, the resulting poset is indeed a
distributive lattice. This happens, for instance, for Dyck,
Motzkin and Schr\"oder paths \cite{FP2}. In all cases, the rank of
a path inside the distributive lattice it lives in is related to
its area. In particular, according to \cite{FP2,FM}, the rank of a
Dyck path of semilength $n$ in his lattice is given by the area
enclosed between the path itself and the path $(UD)^n$ divided by
2. This leads us to the following formulation of the above
conjecture.

\bigskip

\textbf{(A third form of the same) Conjecture.} The distributive
lattice $\mathcal{D}_n$ of Dyck paths of semilength $n$ is
rank-unimodal, for all $n$.

\bigskip

This formulation can be found, for instance, in \cite{FM}, where
the conjecture has been extended also to the Motzkin case (there
is computational evidence supporting this extension). We wish to
point out that there is a case in which rank-unimodality has been
proved, namely that of \emph{Grand Dyck paths} (which are like
Dyck paths, except for the fact that they are allowed to go below
the $x$-axis). Indeed, this is a special case of the unimodality
of Gaussian coefficients (counting integer partitions fitting
inside a rectangle with respect to their size), which has been
proved in several sources using several different methods
\cite{OH,P,Stanl1}. A similar problem concerning compositions
inside a rectangle has been considered in \cite{Sa2}.

In the present paper we have not been able to solve the above mentioned
unimodality conjectures; instead we propose a possible approach
which we haven't been able to find elsewhere, nevertheless we are
strongly convinced that it could prove useful in tackling these
problems.

We start by recalling a particular construction of Dyck paths,
falling into the framework of the so-called \emph{ECO
methodology}. We then show how such a construction suggests a
possible way of decomposing Dyck lattices into saturated chains.
This decomposition gives in turn some hints on what to do to prove
rank-unimodality.

\section{An ECO construction of Dyck paths}

\begin{figure}[ht]
\begin{center}
\begin{tikzpicture}
\begin{scope}[scale=0.2]
\draw (0,0) -- (2,2);
\draw (2,2) -- (3,1);
\draw (3,1) -- (5,3);
\draw (5,3) -- (8,0);
\draw (0,0) [fill] circle (.2);
\draw (1,1) [fill] circle (.2);
\draw (2,2) [fill] circle (.2);
\draw (3,1) [fill] circle (.2);
\draw (4,2) [fill] circle (.2);
\draw (5,3) [fill] circle (.2);
\draw (6,2) [fill] circle (.2);
\draw (7,1) [fill] circle (.2);
\draw (8,0) [fill] circle (.2);

\draw (9,4) -- (17,8);
\draw (9,3) -- (17,4);
\draw (9,2) -- (17,1);
\draw (9,1) -- (17,-3);

\draw (18,9) -- (20,11);
\draw (20,11) -- (21,10);
\draw (21,10) -- (23,12);
\draw (23,12) -- (26,9);
\draw[dotted] (26,9) -- (27,10);
\draw[dotted] (27,10) -- (28,9);
\draw (18,9) [fill] circle (.2);
\draw (19,10) [fill] circle (.2);
\draw (20,11) [fill] circle (.2);
\draw (21,10) [fill] circle (.2);
\draw (22,11) [fill] circle (.2);
\draw (23,12) [fill] circle (.2);
\draw (24,11) [fill] circle (.2);
\draw (25,10) [fill] circle (.2);
\draw (26,9) [fill] circle (.2);

\draw (18,4) -- (20,6);
\draw (20,6) -- (21,5);
\draw (21,5) -- (23,7);
\draw (23,7) -- (25,5);
\draw (27,5) -- (28,4);
\draw[dotted] (25,5) -- (26,6);
\draw[dotted] (26,6) -- (27,5);
\draw (18,4) [fill] circle (.2);
\draw (19,5) [fill] circle (.2);
\draw (20,6) [fill] circle (.2);
\draw (21,5) [fill] circle (.2);
\draw (22,6) [fill] circle (.2);
\draw (23,7) [fill] circle (.2);
\draw (24,6) [fill] circle (.2);
\draw (25,5) [fill] circle (.2);
\draw (28,4) [fill] circle (.2);

\draw (18,0) -- (20,2);
\draw (20,2) -- (21,1);
\draw (21,1) -- (23,3);
\draw (23,3) -- (24,2);
\draw (26,2) -- (28,0);
\draw[dotted] (24,2) -- (25,3);
\draw[dotted] (25,3) -- (26,2);
\draw (18,0) [fill] circle (.2);
\draw (19,1) [fill] circle (.2);
\draw (20,2) [fill] circle (.2);
\draw (21,1) [fill] circle (.2);
\draw (22,2) [fill] circle (.2);
\draw (23,3) [fill] circle (.2);
\draw (24,2) [fill] circle (.2);
\draw (27,1) [fill] circle (.2);
\draw (28,0) [fill] circle (.2);

\draw (18,-5) -- (20,-3);
\draw (20,-3) -- (21,-4);
\draw (21,-4) -- (23,-2);
\draw (25,-2) -- (28,-5);
\draw[dotted] (23,-2) -- (24,-1);
\draw[dotted] (24,-1) -- (25,-2);
\draw (18,-5) [fill] circle (.2);
\draw (19,-4) [fill] circle (.2);
\draw (20,-3) [fill] circle (.2);
\draw (21,-4) [fill] circle (.2);
\draw (22,-3) [fill] circle (.2);
\draw (23,-2) [fill] circle (.2);
\draw (26,-3) [fill] circle (.2);
\draw (27,-4) [fill] circle (.2);
\draw (28,-5) [fill] circle (.2);

\end{scope}
\end{tikzpicture}
\end{center}
\caption{The ECO construction of Dyck paths.\label{ECO}}
\end{figure}
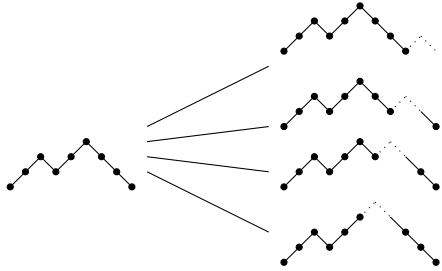

Let $P$ be a Dyck path of length $2n$, and suppose that the length
of its last descent (i.e., of its last sequence of fall steps) is
$k$ (the length of the last descent being simply the number of
fall steps belonging to such a descent). Then we construct $k+1$
Dyck paths of length $2n+2$ starting from $P$ (they will be called
the \emph{sons} of $P$) simply by inserting a peak (that is a rise
step followed by a fall step) in every point of its last descent.
In figure \ref{ECO} it is shown how this construction works.

If one performs such a construction on all Dyck paths of length
$2n$, then it is not difficult to show that every Dyck path of
length $2n+2$ is obtained exactly once.

The above described construction of Dyck paths is well known, and
falls into the framework of the so-called \emph{ECO method} (a
detailed description of which can be found, for instance, in
\cite{BDLPP}).

\section{A decomposition of Dyck lattices into saturated chains}

The above described construction provides a partition of the class
$\mathcal{D}_n$ of Dyck paths of semilength $n$, for all $n$. From
an order-theoretic point of view, it is a partition of
$\mathcal{D}_n$ into saturated chains. We would like to employ
this decomposition to tackle the problem of unimodality in Dyck
lattices. Such a decomposition of the Dyck lattice $\mathcal{D}_n$
will be called the \emph{ECO decomposition} of $\mathcal{D}_n$.

Before starting, we observe that the ECO decomposition of a Dyck lattice
does not possess any of the nice properties usually needed in proving
unimodality: it is not a symmetric chain decomposition, and it is not
even nested (for the notion of a nested chain decomposition, see
for example \cite{G1,G2}).

\bigskip

For any fixed $n$, denote with $P_n$ the $\left( 1+{n-1\choose
2}\right) \times (n-1)$ matrix whose entry $(j,k)$ is the number
of saturated chains of cardinality $k+2$ starting at rank $j$ in
the ECO decomposition of $\mathcal{D}_n$. So, for small values of
$n$, we have the following matrices:
\begin{displaymath}
P_2 =(1),\quad P_3 =\left(%
\begin{array}{cc}
  1 & 0 \\
  0 & 1 \\
\end{array}%
\right) ,\quad P_4 =\left(%
\begin{array}{ccc}
  1 & 0 & 0 \\
  1 & 1 & 0 \\
  0 & 1 & 0 \\
  0 & 0 & 1 \\
\end{array}%
\right) ,\quad P_5 =\left(%
\begin{array}{cccc}
  1 & 0 & 0 & 0 \\
  2 & 1 & 0 & 0 \\
  1 & 2 & 0 & 0 \\
  1 & 1 & 1 & 0 \\
  0 & 1 & 1 & 0 \\
  0 & 0 & 1 & 0 \\
  0 & 0 & 0 & 1 \\
\end{array}%
\right) .
\end{displaymath}

In the sequel we will denote with $P_n ^{(k)}$ the columns of
$P_n$, for $k=0,\ldots ,n-2$.


In a completely analogous way, we also define the $\left(
1+{n\choose 2}\right) \times (n-1)$ matrices $A_n$ whose entry
$(j,k)$ is the number of saturated chains of cardinality $k+2$
ending at rank $j$ in the ECO decomposition of $\mathcal{D}_n$.
Also in this case, $A_n ^{(k)}$ will denote the columns of $A_n$,
for $k=0,\ldots ,n-2$.

Now let $P_n ^{(k)}(x)$ and $A_n ^{(k)}(x)$ be the polynomials
associated with the vectors $P_n ^{(k)}$ and $A_n ^{(k)}$,
respectively; moreover, set
\begin{displaymath}
P_n (x)=\sum_{k=0}^{n-2}P_n ^{(k)}(x),\qquad A_n
(x)=\sum_{k=0}^{n-2}A_n ^{(k)}(x).
\end{displaymath}

So the coefficient of $x^i$ in $P_n (x)$ is the total number of
saturated chains starting at rank $i$, whereas the coefficient of
$x^i$ in $A_n (x)$ is the total number of saturated chains ending
at rank $i$.

Now, in order to link the polynomials $P_n (x)$ and $A_n (x)$ to
our unimodality conjecture, we need to introduce a few more
notations.

Let $r_n ^{(k)}$ denote the number of elements having rank $k$ in
$\mathcal{D}_n$. The polynomial $r_n (x)=\sum_{k}r_n ^{(k)}x^k$ is
called the \emph{rank polynomial} of $\mathcal{D}_n$. We now
introduce the polynomials $s_n (x)$ as follows:
\begin{displaymath}
s_n (x)=\sum_{k}s_n ^{(k)}x^k =\sum_{k}(r_n ^{(k-1)}-r_n
^{(k)})x^k ,
\end{displaymath}
where, by definition, $r_n ^{(-1)}=0$, for all $n$. Obviously, the
unimodality conjecture is equivalent to the following:

\bigskip

\textbf{(Again the same) Conjecture.} For every $n$, there exists
$\bar{k}$ such that $s_n ^{(k)}\leq 0$ for all $k\leq \bar{k}$ and
$s_n ^{(k)}\geq 0$ for all $k>\bar{k}$.

\bigskip

In view of this fact, a possible approach to our problem (which
actually holds in general) consists then of the investigation of
the sign of the coefficients of $s_n (x)$. In this respect, the
following result seems to be of some interest.

\begin{prop}
For all $\ninN$, $s_n (x)=xA_n (x)-P_n (x)$.
\end{prop}

\emph{Proof.}\quad The set of elements having ranks $k-1$ and $k$
in $\mathcal{D}_n$ can be partitioned into three sets (some of
which could be empty): the set of those elements belonging to some
saturated chain crossing both ranks, the set of the top elements
of all chains ending at rank $k-1$ and the set of the bottom
elements starting at rank $k$. The coefficient of $x^k$ in $s_n
(x)$ is clearly given by the difference between the cardinalities
of the second and the third sets described above, whence the
proposition immediately follows.\cvd

There is also a link between the polynomials $P_n ^{(k)}(x)$ and
$A_n ^{(k)}(x)$ which is easy to show and is recorded in the next
proposition.

\begin{prop} For all $\nkinN$, $A_n ^{(k)}(x)=x^{k+1}P_n
^{(k)}(x)$.
\end{prop}

\emph{Proof.}\quad The coefficient of $x^i$ in $x^{k+1}P_n
^{(k)}(x)$ is the number of saturated chains of cardinality $k+2$
starting at rank $i-k-1$. This is clearly the same as the number
of saturated chains of cardinality $k+2$ ending at rank $i$, which
is the coefficient of $x^i$ in $A_n ^{(k)}(x)$.\cvd

As a consequence of the above propositions, we have the following
expression for $s_n (x)$.

\begin{cor}\label{generale} For all $\ninN$, $$s_n (x)=\sum_{k=0}^{n-2}(x^{k+2}-1)P_n ^{(k)}(x).$$
\end{cor}

Everything has been said until this point is valid for \emph{any}
partition into saturated chains of $\mathcal{D}_n$ (and in fact of
any ranked poset). So corollary \ref{generale}, together with a
deep knowledge of the polynomials $P_n ^{(k)}(x)$, could be
helpful in dealing with unimodality.

\bigskip

From now on, we will suppose to work with the ECO decomposition of
$\mathcal{D}_n$ described in the previous section.

\bigskip

We can prove an interesting recurrence for the polynomials $P_n
^{(k)}(x)$, which allows us to interpret the $P_n ^{(k)}(x)$'s as
a polynomial analog of ballot numbers.

\begin{prop}\label{q-ballot} For all $\nkinN$, $$P_n
^{(k)}(x)=x^{k}\cdot (P_{n-1}^{(k-1)}(x)+\cdots
+P_{n-1}^{(n-3)}(x)).$$
\end{prop}

\emph{Proof.}\quad Recall that, for any $i$, the coefficient of
$x^i$ in $P_n ^{(k)}(x)$ is the number of saturated chains having
cardinality $k+2$ and starting at rank $i$ in $\mathcal{D}_n$.
Each of these saturated chains can be uniquely represented by its
minimum, which is a Dyck path $P$ of semilength $n$ ending with
the sequence of steps $UD^{k+1}UD$ (that is, a peak at level 0
preceded by a sequence of $k+1$ consecutive down steps, which is
in turn preceded by an up step). Moreover, as it is shown in
\cite{FP2}, the area $\alpha (P)$ of $P$ is given by
$\alpha (P)=2i+n$. If we remove the last peak, we obtain a
bijection with the set of Dyck paths of semilength $n-1$ ending
with a sequence of precisely $k+1$ $D$ steps and having area
$2i+n-1$. Each of these paths belongs to a different saturated
chains (since they all have the same rank). For any such path $Q$,
the minimum of the saturated chain $Q$ belongs to can be obtained
by simply replacing the last $k+2$ steps of $Q$ (i.e. the sequence
of steps $UD^{k+1}$) with the sequence of steps $D^k UD$. Observe
that, performing such an operation, we are left with a path $R$ in
$\mathcal{D}_{n-1}$ of area $\alpha (R)=\alpha (Q) -2k=2i+n-1-2k$.
This implies that the rank of $R$ in $\mathcal{D}_{n-1}$ is given
by $$r(R)=\frac{\alpha
(R)-(n+1)}{2}=\frac{2i+n-1-n+1-2k}{2}=i-k.$$

Thus we can conclude that the total number of saturated chains
having cardinality $k+2$ and starting at rank $i$ in
$\mathcal{D}_n$ equals the number of saturated chains having
cardinality at least $k+1$ (since its minimum has at least $k$ $D$
steps before the final peak) and starting at rank $i-k$ in
$\mathcal{D}_{n-1}$, whence the thesis follows.\cvd

As an immediate corollary, we also have the following recursion.

\begin{cor} For all
$\nkinN$, with $k\neq 0$, $$P_n ^{(k)}(x)=x\cdot (P_n
^{(k-1)}(x)-x^{k-1}P_{n-1}^{(k-2)}(x)).$$
\end{cor}

\emph{Proof.}\quad A direct application of the above proposition
yields:
\begin{eqnarray*}
P_n ^{(k)}(x)-xP_n ^{(k-1)}(x)&=&x^{k}\cdot
\sum_{i=k-1}^{n-3}P_{n-1}^{(i)}(x)-x^{k}\cdot
\sum_{i=k-2}^{n-3}P_{n-1}^{(i)}(x)\\
&=&-x^k P_{n-1}^{(k-1)}(x).
\end{eqnarray*}
\vspace{-1.2cm}

\cvd

\bigskip

The polynomials $P_n ^{(k)}(x)$'s are not new. They have been
first studied by Carlitz and others \cite{C,CR}, and subsequently
considered also by Krattenthaler \cite{K}. They also found
recursions which look similar to the ones shown in the two above
propositions, however our combinatorial setting is slightly
different. In fact, in the above cited works the combinatorial
meaning of the $P_n ^{(k)}(x)$'s is somehow related to the
distribution of Dyck paths of semilength $n$ with respect to the
area. It is however easy to relate the two approaches. Indeed,
referring to the proof of proposition \ref{q-ballot}, there is a
bijection mapping a saturated chain of the ECO decomposition of
$\mathcal{D}_n$ starting at rank $j$ into a path of
$\mathcal{D}_{n-1}$ at rank $j$. This enables us to interpret the
polynomials $P_n ^{(k)}(x)$'s as describing the distribution of
Dyck paths in $\mathcal{D}_{n-1}$ with respect to rank (i.e. area)
and length of the final descent (i.e. sequence of down steps).

%
%
%

\bigskip

We propose here something new concerning these polynomials, namely
we describe the recursion given by proposition \ref{q-ballot}
using a succession rule with two labels, in the spirit of the ECO
methodology.

\bigskip

In the vector space of polynomials in two indeterminates over the
reals, to be denoted $\mathbf{R}[x,t]$, define the linear operator
$L$ as follows on the canonical basis $(x^\alpha t^\beta )_{\alpha
,\beta \in \mathbf{N}}$:
\begin{align*}
L&:\mathbf{R}[x,t]\longrightarrow \mathbf{R}[x,t]\\
&:x^\alpha t^\beta \longmapsto x^\alpha \cdot \sum_{i=0}^{\beta
+1}x^i t^i .
\end{align*}

We start by noticing an algebraic property of $L$ that will be
useful in what follows.

\begin{lemma}\label{module} $L$ is a homomorphism of the $\mathbf{R}[x]$-module
of polynomials $\mathbf{R}[x][t]$.
\end{lemma}

\emph{Proof.}\quad This is immediate, since
$$
L(x^\alpha t^\beta )=x^\alpha \cdot \sum_{i=0}^{\beta +1}x^i
t^i=x^\alpha L(t^\beta ).
$$
\cvd

%

Now define $P_n (x,t)=L^{n-2}(1)$ (here $L^h$ denotes the
composition of $L$ with itself $h$ times). The following
proposition justifies the introduction of the operator $L$ in our
context.

\begin{prop} For all $n\geq 2$,
\begin{equation}\label{operatore}
P_n (x,t)=\sum_{k=0}^{n-2}P_n ^{(k)}(x)t^k .
\end{equation}
\end{prop}

\emph{Proof.}\quad From the definition of $P_n (x,t)$ it obviously
follows that $L(P_n (x,t))=P_{n+1}(x,t)$. Also, observe that the
degree of $P_n (x,t)$ with respect to $t$ is $n-2$. Set $P_n
(x,t)=\sum_{k=0}^{n-2}Q_n ^{(k)}(x)t^k$, also thanks to the above
lemma, we then obtain:
\begin{align*}
\sum_{k=0}^{n-1}Q_{n+1}^{(k)}(x)t^k &=L\left( \sum_{k=0}^{n-2}Q_n
^{(k)}(x)t^k \right) \\
&=\sum_{k=0}^{n-2}Q_n ^{(k)}(x)L(t^k ) \\
&=\sum_{k=0}^{n-2}Q_n ^{(k)}(x)\left( \sum_{i=0}^{k+1}x^i t^i
\right) \\
&=\sum_{i=0}^{n-2}Q_n ^{(i)}(x)+\sum_{k=1}^{n-1}\left(
\sum_{i=k-1}^{n-2}Q_n ^{(i)}(x)\right) x^k t^k .
\end{align*}

Thus the polynomial sequence $(Q_n ^{(k)}(x))_n$ satisfies the
recursion of proposition \ref{q-ballot}, and of course $Q_2
^{(0)}(x)=P_2 ^{(0)}(x)=1$, whence $Q_n ^{(k)}(x)=P_n ^{(k)}(x)$
which gives the thesis.\cvd

So the operator $L$ encodes the ballot-like recursive generation
of $P_n ^{(k)}(x)$. In the language of the ECO method, $L$ is a
\emph{rule operator} \cite{FP1}. The succession rule described by $L$ (which
turns out to be a two-labelled one) is easily seen to be the
following:
\begin{eqnarray}\label{rule}
\Omega : \left\{ \begin{array}{ll} (0_0 )
\\ (\alpha_{\beta})\rightsquigarrow (\alpha_0 )((\alpha +1)_1 )\cdots ((\alpha +\beta )_{\beta})((\alpha +\beta +1)_{\beta +1})
\end{array}\right. .
\end{eqnarray}

\begin{figure}[ht]
\begin{center}
\begin{tikzpicture}
\begin{scope}[scale=1]
\draw (0,0) node{$0_0$};
\draw (-2,-1) node{$0_0$};
\draw (2,-1) node{$1_1$};
\draw (-3,-2) node{$0_0$};
\draw (-1,-2) node{$1_1$};
\draw (1,-2) node{$1_0$};
\draw (2,-2) node{$2_1$};
\draw (3,-2) node{$3_2$};
\draw (-4,-3) node{$0_0$};
\draw (-3,-3) node{$1_1$};
\draw (-2,-3) node{$1_0$};
\draw (-1.5,-3) node{$2_1$};
\draw (-1,-3) node{$3_2$};
\draw (0.5,-3) node{$1_0$};
\draw (1,-3) node{$2_1$};
\draw (1.5,-3) node{$2_0$};
\draw (2,-3) node{$3_1$};
\draw (2.5,-3) node{$4_2$};
\draw (3,-3) node{$3_0$};
\draw (3.5,-3) node{$4_1$};
\draw (4,-3) node{$5_2$};
\draw (4.5,-3) node{$6_3$};
\draw (-0.2,-0.2) -- (-1.8,-0.8);
\draw (0.2,-0.2) -- (1.8,-0.8);
\draw (-2.2,-1.2) -- (-2.8,-1.8);
\draw (-1.8,-1.2) -- (-1.2,-1.8);
\draw (1.8,-1.2) -- (1.2,-1.8);
\draw (2,-1.2) -- (2,-1.8);
\draw (2.2,-1.2) -- (2.8,-1.8);
\draw (-3.2,-2.2) -- (-3.8,-2.8);
\draw (-3,-2.2) -- (-3,-2.8);
\draw (-1.2,-2.2) -- (-1.8,-2.8);
\draw (-1.1,-2.2) -- (-1.5,-2.8);
\draw (-1,-2.2) -- (-1,-2.8);
\draw (0.9,-2.2) -- (0.6,-2.8);
\draw (1,-2.2) -- (1,-2.8);
\draw (1.9,-2.2) -- (1.4,-2.8);
\draw (2,-2.2) -- (2,-2.8);
\draw (2.1,-2.2) -- (2.4,-2.8);
\draw (-1.2,-2.2) -- (-1.8,-2.8);
\draw (3,-2.2) -- (3,-2.8);
\draw (3.1,-2.2) -- (3.4,-2.8);
\draw (3.2,-2.2) -- (3.9,-2.8);
\draw (3.3,-2.2) -- (4.4,-2.8);
\end{scope}
\end{tikzpicture}
\end{center}
\caption{The generating tree associated with $\Omega$.\label{gentree}}
\end{figure}
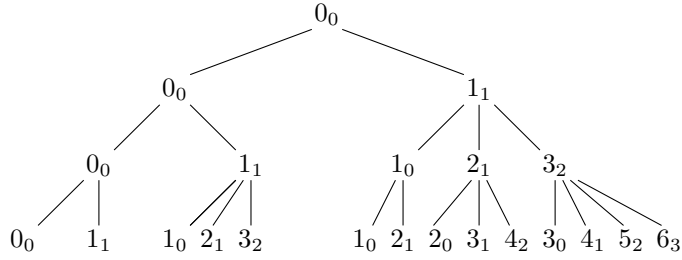

The first levels of the generating tree of this rule are depicted
in figure \ref{gentree}.

The infinite matrix describing the distribution of the labels at the various levels
of the generating tree is called the \emph{ECO matrix} of the succession rule in \cite{DFR}.
In our case, for the rule in (\ref{rule}), the first lines of such a matrix are the following:

\bigskip

\scriptsize
\begin{center}
\begin{tabular}{|c|c|cc|cc|ccc|ccc|ccc|cc|}
\hline

\tiny{labels}$\rightarrow$& $0_0$ & $1_0$ & $1_1$ & $2_0$ & $2_1$ &
$3_0$ & $3_1$ & $3_2$ & $4_0$ & $4_1$ & $4_2$ & $5_0$ & $5_1$ & $5_2$ & $6_0$ & $\ldots$\\
\tiny{at level} $\downarrow$ &&&&&&&&&&&&&&&&\\

\hline

&&&&&&&&&&&&&&&&\\
0&1&&&&&&&&&&&&&&&\\
1&1&0&1&&&&&&&&&&&&&\\
2&1&1&1&0&1&0&0&1&&&&&&&&\\
3&1&2&1&1&2&1&1&1&0&1&1&0&0&1&0&$\cdots$\\
4&1&3&1&3&3&3&3&1&2&3&2&1&2&2&1&$\cdots$\\
5&1&4&1&6&4&7&6&1&7&7&3&5&7&4&5&$\cdots$\\
6&1&5&1&10&5&14&10&1&17&14&4&16&17&7&16&$\cdots$\\
7&1&6&1&15&6&25&15&1&35&25&5&40&35&11&43&$\cdots$\\
$\vdots$&$\vdots$&$\vdots$&$\vdots$&$\vdots$&$\vdots$&$\vdots$&$\vdots$&$\vdots$&$\vdots$&
$\vdots$&$\vdots$&$\vdots$&$\vdots$&$\vdots$&$\vdots$&$\vdots$\\
\hline
\end{tabular}
\end{center}
\normalsize

\bigskip

We can also find a recursion for the columns of the above ECO matrix.
Since it directly depends on the succession rule (\ref{rule}), the proof is left to the reader.

\begin{prop} Denote with $\mathcal{C}_{k_i}(x)$ the generating function of the column
associated with the label $k_i$. Then the following recursion holds:
$$
\mathcal{C}_{k_i}(x)=x\cdot \sum_{j\geq i-1}\mathcal{C}_{(k-i)_j}(x).
$$
\end{prop}

\section{Conclusions and further work}

We have addressed the problem of the rank-unimodality of Dyck
lattices and, as announced at the beginning, we have not been able
to solve it. Nevertheless, we hope to have provided some
interesting insight to the problem, as well as an original way to
tackle it. We hope that someone more skillful than us will be able
to further develop these ideas to eventually find the desired
unimodality proof (maybe by finding more structural and
enumerative properties of the ECO decomposition of Dyck lattices).

\bigskip

We remark that a similar approach can be considered, in which we make use of a
different decomposition of Dyck lattices. More precisely, for any Dyck lattice
$\mathcal{D}_n$, one can consider the sublattices $D_{n,k}$ consisting of
all paths starting with exactly $k$ up steps (for $1\leq k\leq n$).
This is not of course a decomposition into chains. What is interesting about
this decomposition is that, for any $k$, $\mathcal{D}_{n,k}$ has a natural embedding
into $\mathcal{D}_{n,k-1}$ (see \cite{FP2}), a fact that could be useful in proving
unimodality.

\bigskip

We also notice that, in order to prove unimodality for an integer sequence,
an extremely useful information is the position of the maximum. In our case,
we even do not know where the maximums are located, and how they depend on $n$.
This is a related open problem that deserves to be solved.

\bigskip

We close this paper by observing that the approach we have presented here to study
unimodality of Dyck lattices is suggested by the ECO construction of Dyck paths,
and that, in the end, to prove rank-unimodality of Dyck lattices it would be enough
to prove that, at each level of the associated generating tree,
the distribution of the labels $k$ of the succession rule (\ref{rule}) is unimodal
(where $k=\sum_i k_i$). We are thus led to formulate a \emph{unimodality problem} for
succession rules in general: given a succession rule, is it possible to find some
(necessary and/or sufficient) conditions on the rule itself for the unimodality of the
integer sequences describing the distribution of the labels
at the various levels of the associated generating tree?
This is a problem of independent interest, which seems not to have been previously considered.

To give just the flavor of this kind of investigations, we prove a
simple result concerning finite succession rules. We will say that
a succession rule is \emph{unimodal} when the integer sequences
appearing in the rows of the associated ECO matrix are all
unimodal. Moreover, a succession rule deprived of its axiom will
be called a \emph{family of succession rules} (since each choice
of the axiom determines a distinct succession rule).

\begin{teor} Fix $\ninN$ and let
$$
\Omega :(k)\rightsquigarrow (1)^{a_1 ^{(k)}}(2)^{a_2 ^{(k)}}\cdots (n)^{a_n ^{(k)}}
$$
be a family of succession rules. For any $b\in \mathbf{N}$, denote with $\Omega^{(b)}$
the succession rule of the family $\Omega$ having axiom $(b)$. If
$\Omega^{(b)}$ is unimodal for all $b\leq n$, then the sequences
$(a_1^{(k)},a_2^{(k)},\ldots ,a_n^{(k)})$ are unimodal for all $k\leq n$.
Vice versa, if, for all $k\leq n$, the sequences $(a_1^{(k)},a_2^{(k)},\ldots ,a_n^{(k)})$
are unimodal and all have the maximum in the same position,
then $\Omega^{(b)}$ is unimodal for all $b\leq n$.
\end{teor}

\emph{Proof.}\quad Suppose that all the rules of the family $\Omega$ are unimodal and, for any given $k\leq n$,
consider the succession rule $\Omega ^{(k)}$:
the sequence $(a_1^{(k)},a_2^{(k)},\ldots ,a_n^{(k)})$ appears in the first row of its
ECO matrix, so such a sequence is unimodal (since $\Omega ^{(k)}$ is unimodal by hypothesis).
On the other hand, suppose that, for all $k\leq n$,
the sequences $(a_1^{(k)},a_2^{(k)},\ldots ,a_n^{(k)})$ are unimodal and also that
all have maximum at index $i$. Fixed $b\leq n$, let $\beta_1 ,\beta_2 ,\ldots ,\beta_n$
be the $l$-th row of the ECO matrix of $\Omega_b$. The recursion determined by $\Omega_b$
allows to express any row of the ECO matrix in terms of the previous one.
Therefore, if $\alpha_1 ,\alpha_2 ,\ldots ,\alpha_n$ is the $l-1$-th row, we have:
\begin{align*}
\beta_1 &= \alpha_1 a_1 ^{(1)}+\cdots +\alpha_n a_1 ^{(n)}\\
&\vdots\\
\beta_n &= \alpha_1 a_n ^{(1)}+\cdots +\alpha_n a_n ^{(n)}
\end{align*}

The unimodality hypotheses on the sequences $(a_1^{(k)},a_2^{(k)},\ldots ,a_n^{(k)})$
immediately implies that $\beta_1 \leq \cdots \leq \beta_{i-1}\leq \beta_i \geq \beta_{i+1}
\geq \cdots \geq \beta_n$.\cvd

\emph{Remark.}\quad We wish to point out that, in the above theorem,
the hypothesis concerning the position of the maximum is essential.
Indeed, consider the following family of finite succession rules:
$$
\Omega :\left\{ \begin{array}{llll}
(1)\rightsquigarrow (3)\\
(2)\rightsquigarrow (2)(3)\\
(3)\rightsquigarrow (3)(4)(4)\\
(4)\rightsquigarrow (1)(1)(2)(2)
\end{array}\right. .
$$

One immediately notices that, for each $k\leq 4$, the sequence of the productions of
$(k)$ is unimodal. However, it is not difficult to realize that none of the
succession rules $\Omega^{(b)}$ is unimodal (for $b\leq 4$).

\end{document}